\documentclass[11pt]{article}
\usepackage{amsmath}
\usepackage{amssymb}
\usepackage{mathptmx}
\usepackage{amsfonts}
\usepackage[mathscr]{eucal}
\usepackage[dvips]{graphicx}
\usepackage{color}

\newcommand{\mbf}[1]{\ensuremath{\mathbf{#1}}}
\newcommand{\ms}[1]{\ensuremath{\mathscr{#1}}}
\newcommand{\sqz}[1]{\ensuremath{d\Gamma({#1}) }}
\newcommand{\sqzb}[1]{\ensuremath{
d\Gamma_{\textrm{b}}({#1}) }}

\newcommand{\Hb}{H_{\textrm{b}}}
\newcommand{\Hpar}{H_{\textrm{p}}}
\newcommand{\HI}{H_{\textrm{I}}}

\newcommand{\effect}{\textrm{eff}}
\newcommand{\kuuhaku}{\; \;}

\newcommand{\tens}{\otimes}

\newcommand{\nstens}{\otimes^{n}_{\textrm{s}}}
\newcommand{\eltwo}{L^{2}(\mathbf{R}^{3})}
\newcommand{\eltwod}{L^{2}(\mathbf{R}^{d})}

\newcommand{\eltwodNx}{L^{2}(\mathbf{R}^{dN}_{\mbf{x}})}
\newcommand{\eltwodk}{L^{2}(\mathbf{R}^{d}_{\mbf{k}})}

\newcommand{\Rd}{\mathbf{R}^{d} }
\newcommand{\RdN}{\mathbf{R}^{dN} }

\newcommand{\Fb}{\mathscr{F}_{\textrm{b}} ( L^{2}(\mathbf{R}^{d}) )  }

\newcommand{\Omegab}{ \Omega_{\textrm{b}} }
\newcommand{\POmegab}{ P_{\Omega_{\textrm{b}} }}
\newcommand{\bos}{\textrm{b}}
\newcommand{\fin}{\textrm{fin}}
\newcommand{\eff}{\textrm{eff}}

\newcommand{\phattwo}[2]{\ensuremath{\hat{p}^{#2}_{#1}}}
\newcommand{\phatone}[1]{\ensuremath{\hat{p}^{#1} }}

\newtheorem{theorem}{Theorem}[section]
\newtheorem{proposition}[theorem]{Proposition}
\newtheorem{lemma}[theorem]{Lemma}
\newtheorem{corollary}[theorem]{Corollary}
\newtheorem{remark}{Remark}[section]

\begin{document}

\title{Scaling Limits for the System of Semi-Relativistic Particles \\ 
  Coupled to  a Scalar Bose Field}

\author{Toshimitsu TAKAESU }
\date{ }
\maketitle

\begin{center}
\textit{Faculty of Mathematics, Kyushu University,\\  Fukuoka, 812-8581, Japan }
\end{center}

\begin{quote}
\textbf{Abstract}. In this paper the Hamiltonian for the system of  semi-relativistic particles interacting with a scalar bose field  is investigated. A scaled total Hamiltonian of the system
 is defined and its scaling limit is considered. Then  the semi-relativistic Schr\"odinger 
  operator with an effective potential is derived.\\
 
{ \small
key words : Spectral analysis, Relativistic Schr\"odinger operator, Quantum field theory.}
\end{quote}

\section{Introduction}
In this paper we consider the Hamiltonian of the system of $N$ particles linearly coupled to a  scalar bose field. 
We assume that particles obey the semi-relativistic Schr\"{o}dinger operator
\[
    \Hpar \; \; = \; \;    \sum_{j=1}^{N} \sqrt{-\triangle_{j} + M^2 } ,
\]
 where $M>0 $ is a rest mass. There has been   many results on the  spectral properties of $\Hpar$. 
 Refer to  e.g. \cite{We74, He77, Da84, De90, Ge91, Um95}, and see also \cite{SSS10}.
The free Hamiltonian  $\Hb$ of the scalar bose field is defined by the second quantization of  
 the multiplication operator $\omega $, which is formally expressed by
\[
\Hb \; \;  = \; \;  \int_{\Rd} \omega (\mbf{k} ) a^\ast (\mbf{k}) a(\mbf{k}) d\mbf{k}
\] 
The state space of the interacting system is defined by
 $ \ms{H} \; = \; \eltwodNx \tens \ms{F}_{\bos} (\eltwodk ) \; $ where 
 $ \ms{F}_{\bos} (\eltwodk )  $ is the boson Fock space on $\eltwodk$. The total Hamiltonian is 
 given by 
\begin{equation}
H \; \; = \; \;  
 \Hpar \tens I  \;  \;   + \; \; I \tens  H_{\bos} \;   + \;   \kappa \HI  , \qquad \qquad \kappa \in \mbf{R}.
    \end{equation}
Here the interaction $\HI$ is denoted by formally 
\[
  \HI \; \; = \; \; \sum_{j=1}^N \frac{1}{\sqrt{2}}
  \int_{\Rd} \left( \frac{}{} \overline{f_{\mbf{x_j }}(\mbf{k})} \tens a(\mbf{k}) + f_{\mbf{x_j }}(\mbf{k}) \tens a^{\ast} (\mbf{k}) \right) d \mbf{k},
 \]
where 
$ f_{\mbf{x}} $ is an multiplication operator on $L^2 (\Rd_{\mbf{x}})$.

$\quad $\\
We consider the scaled Hamiltonian 
\begin{equation}
H (\Lambda )\; \; = \; \;  
 \Hpar \tens I \;  \;
  + \; \; \Lambda^2 \; I \tens   H_{\bos} \;   + \;   \kappa \Lambda \HI  ,  \qquad \qquad \Lambda > 0 . \label{HLambda}
    \end{equation}
  We investigate the asymptotic behavior of $H (\Lambda)$ as $\Lambda \to \infty$. 
  The unitary evolution $e^{ -itH(\Lambda )}$ generated by $H (\Lambda )$  is given by 
{\large
 \[
 e^{ -itH(\Lambda )} \; \; 
= \; \; 
 e^{\; \; -i \Lambda^2 \, t \left( \; 
 \sum\limits_{j=1}^N \sqrt{ \left( \frac{\hat{\mbf{p}}_{j} }{\Lambda^2} \right)^2  + \left( \frac{M}{\Lambda^2} \right)^2 }
  \;   \;  +   \; \;  \Hb \; \; + \; \;  \left( \frac{\kappa}{\Lambda} \right)   \HI
\; \right) }.
\]
}Here $\Lambda^2 t $  denotes the scaled time, $\; \; {\Lambda^{-2}}\hat{\mbf{p}} \; $ the scaled momentum for $\hat{\mbf{p}} = - i \mbf{\nabla}$,
  $\; \Lambda^{-2} M $ the scaled mass, and $\; \Lambda^{-1} \kappa$ the scaled coupling constant.  As far as we know scaling limits
of
  the 
Hamiltonians of the form
(\ref{HLambda}) is  initiated by E. B. Davies \cite{Da79}, where
$H(\Lambda)$ with semi-relativistic Schr\"odinger operator 
replaced by a standard Schr\"odinger operator
is considered and
a  Sch\"{o}dinger operator with an effective potential
is derived as $\Lambda\to \infty$.
This model is called the Nelson model,
and our result can be regarded as a semi-relativistic version  of \cite{Da79}.
 In \cite{Ar90},
a general theory of
scaling limits is  established and it is applied to scaling limits of a spin-boson model and
 non-relativistic QED models.  In \cite{Hi98}, by removing
ultraviolet cutoffs and taking a scaling limit
of the Nelson model {\it simultaneously},
 a Schr\"odinger operator
with the Yukawa potential
or the Coulomb potential
is derived.
Refer to see also
  \cite{Hi93, Su07-1, Su07-2, Oh09, Ta09}.

$\quad$ \\
In the main theorem, it is shown that for $z \in \mbf{C} \backslash \mbf{R}$,
  \begin{equation}
s-\lim_{\Lambda \to \infty} \left( H (\Lambda ) 
  -z \right)^{-1} \; 
  =  \;  \left(  \sum_{j=1}^{N} \sqrt{-\triangle_{j} + M^2 }
 \; \; + \; \;   V_{\eff} (\mbf{x}_{1} , \cdots , \mbf{x}_{N}) \;  \; -\; z \right)^{-1}  \; \POmegab , 
\end{equation}
where 
\[
V_{\eff} (\mbf{x}_{1} , \cdots , \mbf{x}_{N}) \; \; = \; \; - \frac{\kappa^2}{4} \sum_{j,l } \int_{\Rd} \frac{\overline{f_{\mbf{x}_j}(\mbf{k}) } 
f_{\mbf{x}_{l}} (\mbf{k})  + \overline{f_{\mbf{x}_l}(\mbf{k}) } 
f_{\mbf{x}_{j}} (\mbf{k}) 
}{\omega(\mbf{k})} d \mbf{k} \; ,
\]
and $\POmegab$ is the projection onto the closed subspace spanned by the Fock vacuum $\Omegab$ of the bose field. \\
$\quad$ \\
$\quad$ \\
For the strategy of the proof of the main theorem, we use a unitary transformation, 
called \textit{the dressing transformation}. Then we apply the general theory investigated in \cite{Ar90} to 
 the unitary transformed Hamiltonian $U(\Lambda)^{-1} H (\Lambda) U (\Lambda )$, and the consider the 
 asymptotic behavior of    $U(\Lambda)^{-1} H (\Lambda) U (\Lambda )$ as $\Lambda \;  \to \; \infty$. \\
 $\quad$ \\

$\quad$ \\
This paper is organized as follows. In Section 2, the theory of boson Fock space is described. Then  
the total  state space and the total Hamiltonian is defined, and the main results are stated.  In Section 3, the proof of the main theorem is given.

%%%%%%%%%%%%%%%%%%%%%%%%%%%%%%%%%%%%%%%%%%%%%%%%%%%%%%%%%%%%%%%%%%%%%%%%%%%%%%%%%%%%%%%%%%%%%%%%%%%%%%%%%%%%%%%%%%%%%%%%%%%%%%%%%%%%%%%%%%%%%%%%%%%%%%%%%%%%%%%%%%%%%%%%%%%%%%%%%%%%%%%%%%%%%%%%%%%%%%%%%%%%%%%%%%%%%%%%%%%%%%%%%%%%%%%%%%%%%%%%%%%%%%%%%%%%%%%%%%%%%%%%%%%%%%%%%%%%%%%%%%%%%%%%%%%%%%%%%%%%%%%%%%%%%%%%%%%%%%%%%%%%%%%%%%%%%%%%%%%%%%%%%%%%%%%%%%%%%%%%%%%%%%%%%%%%%%%%%%%%%%%%%%%%%%%%%%%%
\section{Main Results}
\subsection{Boson Fock Spaces}
In this subsection we give the mathematically rigorous definition of the bose field. 
The state space of the bose field is given by the boson Fock space 
$ \Fb \;   = \; \oplus_{n=0}^{\infty}
 (  \nstens \eltwod ) $,
where $\nstens \eltwod$ denotes the $n$-fold symmetric tenser
 product of $\eltwod \; $ with $\tens_{s}^{0} \eltwod := \mbf{C}$.
 The Fock vacuum is defined by $\Omegab = \{ 1, 0, 0, \cdots \} \in \Fb$.
   The finite particle subspace
 $\ms{F}_{\bos}^{\; \fin} (\ms{D})$ on the subspace $\ms{D} \subset \eltwod$ is defined by the set  of $\Psi = \{\Psi^{(n)} \}_{n=0}^{\infty} \; $ satisfying  that
   $\Psi^{(n)}  \in \nstens \ms{D}  $, $ n  \geq 0 $, and
    $\Psi^{(n')}  = 0  $ for all $ n' > N$ with  some $N \geq 0$.
 Let $a( \xi ), \; \xi  \in \eltwo $, and $a^{\ast}(\eta ), \; \eta \in \eltwod $, be the annihilation operator and the creation operator on $\Fb$, respectively.
  Then they satisfy the canonical commutation relations on  $\ms{F}_{\bos}^{\fin}(  \eltwod )$ :
\[
    [ \, a( \xi ), \, a^{\ast} (\eta )  ] = (\xi , \eta  ), \qquad
    [ \, a  (\xi ), \, a (\eta  )  ] = [ a^{\ast}(\xi  ), \,  a^{\ast}(\eta ) ] =0 .
\]
Let $S $ be a  self-adjoint operator on $\eltwod$. The  second quantization of $S$ is defined by
\[
 \sqz{S}  \; = \bigoplus_{n=0}^{\infty} \left(
\sum_{j=1}^{n} ( I \tens \cdots I \tens \underbrace{S}_{jth} \tens I   \cdots  \tens I )\right),
\]
 For $\eta \in \ms{D} (S^{-1/2})$, it is seen that $a(\eta) $ and $a^{\ast} (\eta)$ are  relatively bounded with respect to
$\sqz{S}^{1/2} $  with the bound
\begin{align}
& \| a (\eta ) \Psi  \| \leq \| S^{-1/2} \eta  \|  \, \|  \sqz{S}^{1/2} \Psi \| ,
    \qquad \qquad  \qquad  \quad  \Psi  \in  \ms{D} ( \sqz{S}^{1/2} ),    \label{bounda}   \\
& \| a^{\ast}(\eta ) \Psi  \| \leq \| S^{-1/2} \eta \|
\|   \sqz{S}^{1/2} \Psi \| + \| \eta  \| \| \Psi \| ,  \qquad   \Psi  \in  \ms{D} (
 \sqz{S}^{1/2} ) .
\label{boundad}
\end{align}
The field operator and its conjugate operator are defined by 
\[
 \phi (\xi) \; = \; \frac{1}{\sqrt{2}}
 \left( \frac{}{} a (\xi ) \; + \;  a^{\ast} (\xi ) \right) ,  \qquad 
  \Pi (\eta ) \; = \; \frac{i}{\sqrt{2}}
 \left( \frac{}{}   - a (\eta ) \; + \;  a^{\ast} (\eta ) \right) .
\]

%%%%%%%%%%%%%%%%%%%%%%%%%%%%%%%%%%%%%%%%%%%%%%%%%%%%%%%%%%%%%%%%%%%%%%%%%%%%%%%%%%%%%%%%%%%%%%%%%%%%%%%%%%%%%%%%%%%%%%%%%%%%%%%%%%%%%%%%%%%%%%%%%%%%%%%%%%%%%%%%%%%%%%%%%%%%%%%%%%%%%%%%%%%%%%%%%%%%%%%%%%%%%%%%%%%%%%%%%%%%%%%%%%%%%%%%%%%%%%%%%%%%%%%%%%%%%%%%%%%%%%%%%%%%%%%%%%%%%%%%%%%%%%%%%%%%%%%%%
\subsection{Main Theorem}
In this subsection we define  the total Hamiltonian and state the main results.
The state space of the system for the $N$-particles coupled to bose field is defined by 
\[
\ms{H} \; = \;  \eltwodNx \tens \ms{F}_{\bos} (\eltwodk ) . 
\]
The free Hamiltonian of particles and the bose field are defined by 
\[
\Hpar \; \; =  \; \; \;  \sum_{j=1}^{N} \sqrt{-\triangle_{j} + M^2 }  \; ,  
 \qquad \qquad \Hb  \; \;  = \; \; \sqzb{\omega} ,  \qquad \qquad 
\]
where $M> 0 $ is a rest mass and $\omega $ denotes the multiplication operator by the function $\omega (\mbf{k})$, 
which  describes the  energy of the  boson with momentum $\mbf{k}$. 
 We assume 
 the following condition : 
\begin{quote}
\textbf{(A.1)} $\omega $ is non-negative.
\end{quote}
The interaction $\HI $   is defined by
\[
\HI \; \; = \sum_{j=1}^N  \; \phi (f_{\mbf{x}_{j}})  , 
\]
where $f_{\mbf{x}}$ is  the multiplication operator satisfying the following condition : 
\begin{quote} $\; $ \\
\text\bf{(A.2)}
\[
\sup_{\mbf{x} \in \Rd} \; 
\int_{\Rd }   | f_{\mbf{x}} (\mbf{k}) |^2 d \mbf{k} \; <  \; \infty ,  \qquad  \text{and} \qquad 
  \sup_{\mbf{x} \in \Rd} \;
\int_{\Rd } \frac{}{}  \frac{| f_{\mbf{x}} (\mbf{k}) |^2  }{\omega (\mbf{k})}    d \mbf{k} \; <  \; \infty . 
\]
\end{quote}

$\quad$ \\
The total Hamiltonian of this system is given by
\[
H \; \; = \; \;  H_{0} \; \; + \; \; \kappa \; \HI  ,
\]
where $\; H_{0 } \;  = \; \Hpar \tens I  \; +  \; I \tens \Hb $.
By (\ref{bounda}), (\ref{boundad}), and the assumption (\textbf{A.2}), it is seen that
 the $\HI$ is relatively bounded with respect to $I \tens \Hb^{1/2}$. 
Hence $\HI$ is relatively bounded with respect to $H_{0}$ with infinitely small bound. Then the  Kato-Rellich theorem shows that $H$ is self-adjoint and essentially self-adjoint any core of $H_{0}$. 
 Then in particular, $H$ is essentially self-adjoint on $  \ms{D}_{0} \; \; = \; \; C_{0}^{\infty} (\RdN )  \hat{\tens} \ms{F}_{\bos}^{\fin} (\ms{D} (\omega )) $, where $\hat{\tens}$ denotes the algebraic tensor product.

$\quad$ \\
Let us introduce the scaled total Hamiltonian
\[
H (\Lambda ) \; \; = \; \; H_{0}(\Lambda ) \; \; + \; \; \kappa \Lambda \; \HI  ,
\]
where $\; H_{0 } (\Lambda ) \;  = \; \Hpar  \tens I \; +  \Lambda^2 I \tens  \Hb $. 
We introduce an additional assumption on the interaction.\\
\begin{quote}
\textbf{(A.3)} $ \sup\limits_{\mbf{x} \in \Rd} \; 
\int_{\Rd }   \frac{|f_{\mbf{x}} (\mbf{k})|^2}{\omega (\mbf{k})^2} d \mbf{k} \; <  \; \infty $, $\quad $
 and $\quad $ 
$ \; \sup\limits_{\mbf{x} \in \Rd} \; 
\int_{\Rd }   \frac{|f_{\mbf{x}} (\mbf{k})|^2}{\omega (\mbf{k})^3} d \mbf{k} \; <  \; \infty $.
\end{quote}
\begin{quote}
\textbf{(A.4)}  $ \; \sup\limits_{\mbf{x} \in \Rd} \; 
 \int_{\Rd }   \frac{| \partial_{x^\nu} f_{\mbf{x}} (\mbf{k})|^2}{\omega (\mbf{k})^2} d \mbf{k} \; <  \; \infty \; $, 
  $ \; \sup\limits_{\mbf{x} \in \Rd} \; 
 \int_{\Rd }   \frac{| \triangle f_{\mbf{x}} (\mbf{k})|^2}{\omega (\mbf{k})^2} d \mbf{k} \; <  \; \infty \; $ and 
$ \;  ( \frac{\partial_{x^\nu } f_{\mbf{x}}}{\omega} , \frac{f_{\mbf{y}}}{\omega}) \in \mbf{R} , \; \;    \; \;  \mbf{x} , \mbf{y} \in \Rd$.
\end{quote}
$\quad$\\
Under the condition $ \;  ( \frac{\partial_{x^\nu } f_{\mbf{x}}}{\omega} , \frac{f_{\mbf{y}}}{\omega}) \in \mbf{R} $ in \textbf{(A.4)}, it follows that $ [ \Pi ( \frac{\partial_{x^\nu } f_{\mbf{x}}}{\omega} ) , 
 \Pi (\frac{f_{\mbf{y}}}{\omega} )] = 0$, $\mbf{x} , \mbf{y}  \in \Rd$. 
\begin{remark}
 Let us define that  
 $\; f_{\mbf{x}} (\mbf{k})\; = \; \frac{\chi_{ R} (|\mbf{k}|)}{ \sqrt{ \omega
 (\mbf{k})}} \; e^{-i \mbf{k} \cdot \mbf{x}} \; $ with $\omega (\mbf{k} ) = \omega (-\mbf{k}) $. Here   $ \; \chi_{ R} $ denotes the characteristic function on $[0,\; R)$. Then  the conditions \textbf{(A.1)}-\textbf{(A.4)} are satisfied, and the  interaction $\HI$ is formally expressed by
\[
\HI \; \; = \; \; \sum_{j=1}^N \, \int_{\Rd} \; \frac{\chi_{R} (|\mbf{k}|)}{ \sqrt{ 2 \omega (\mbf{k})}} 
 \left( \frac{}{} a(\mbf{k}) e^{i \mbf{k} \cdot \mbf{x}_{j}}  \; + \; a^{\ast}(\mbf{k}) e^{-i \mbf{k} \cdot \mbf{x}_j } \right)
  d \mbf{k} .
\]
\end{remark}
$\quad$ \\
$\quad$ \\
The main theorem in this paper is as follows
\begin{theorem} \label{MainTheorem}
Assume (\textbf{A.1})-(\textbf{A.4}). Then for $z \in \mbf{C} \backslash \mbf{R}$ it follows that 
\[
 s-\lim_{\Lambda \to \infty} \left( H (\Lambda ) 
  -z \right)^{-1} \; 
  =  \;  \left( \Hpar \; \; 
  + \; \;  V_{\eff} (\mbf{x}_{1}, \cdots , \mbf{x}_{n}) \; \;  -\; \; z \right)  \tens \; \POmegab , 
  \]
where 
\[
V_{\eff} (\mbf{x}_{1}, \cdots , \mbf{x}_{n})\; \; = \; \; - \frac{\kappa^2}{4} \sum_{j,l } \int_{\Rd} 
\frac{\overline{f_{\mbf{x}_j}(\mbf{k}) } f_{\mbf{x}_{l}} (\mbf{k}) + 
\overline{f_{\mbf{x}_l}(\mbf{k}) } f_{\mbf{x}_{j}} (\mbf{k}) }{\omega(\mbf{k})}   d \mbf{k},
\]
and $\POmegab$ is the projection onto the closed subspace spanned by the Fock vacuum $\Omegab$. \\
\end{theorem}
\begin{remark}
When $\; f_{\mbf{x}} (\mbf{k})\; = \; \frac{\chi_{ R} (|\mbf{k}|)}{ \sqrt{ \omega
 (\mbf{k})}} \; e^{-i \mbf{k} \cdot \mbf{x}} $ with $\omega (\mbf{k} ) = \omega (-\mbf{k}) $, the effective potential is given by 
\[
V_{\eff} (\mbf{x}_{1}, \cdots , \mbf{x}_{n})\; \; = \; \; - \frac{\kappa^2}{2} \sum_{j,l } \int_{\Rd} \frac{|\chi_{R}(\mbf{k})|^2 }{\omega(\mbf{k})^2} e^{\, -i \mbf{k} \cdot (\mbf{x}_{j} - \mbf{x}_{l}) }  d \mbf{k},
\] 
\end{remark}
By using the  norm convergence theorem considered in (\cite{Su07-2} ; Lemma 2.7), the the next corollary follows.
\begin{corollary}
Assume (\textbf{A.1})-(\textbf{A.4}). Then  it follows that 
\[
 s-\lim_{\Lambda \to \infty} e^{-itH (\Lambda )}
 \left(  I \tens  \POmegab \right) \; 
  = e^{ -it \left( \;  \Hpar 
  \; \; + \; \; V_{\eff}   (\mbf{x}_{1}, \cdots , \mbf{x}_{n}) \right) }  \; \tens \POmegab . 
  \] 
\end{corollary}

%%%%%%%%%%%%%%%%%%%%%%%%%%%%%%%%%%%%%%%%%%%%%%%%%%%%%%%%%%%%%%%%%%%%%%%%%%%%%%%%%%%%%%%%%%%%%%%%%%%%%%%%%%%%%%%%%%%%%%%%%%%%%%%%%%%%%%%%%%%%%%%%%%%%%%%%%%%%%%%%%%%%%%%%%%%%%%%%%%%%%%%%%%%%%%%%%%%%%%%%%%%%%%%%%%%%%%%%%%%%%%%%%%%%%%%%%%%%%%%%%%%%%%%%%%%%%%%%%%%%%%%%%%%%%%%%%%%%%%%%%%%%%%%%%%%%%%%%%
\section{Proof of Main Theorem}
The outline  of the proof of Theorem 2.1 is as follows. 
A unitary transformation  $U(\Lambda)$, called the dressing transformation, is defined and we consider the 
 unitarily transformed Hamiltonian 
 $U(\Lambda )^{-1} H (\Lambda ) U (\Lambda)$. Then we apply the general theory on  scaling limits
  in \cite{Ar90} to $U(\Lambda )^{-1} H (\Lambda ) U (\Lambda)$.

$\quad$ \\ 
Under the condition (\textbf{A.3}),  the following unitary operator can be  defined : 
\[  
U(\Lambda ) \; = \; 
 e^{i  \left( \frac{\kappa}{\Lambda} \right)  \sum\limits_{j=1}^{N} 
 \Pi ( \frac{f_{\mbf{x}_{j} }}{\omega} ) } .
\] 
It is  seen that on the finite particle subspace 
\begin{align}
& [\Pi (\xi ) , \Hb ] \; = \; -i \, \phi (\omega  \xi ) , \qquad \xi \in \ms{D} (\omega) ,  \label{PiHb} \\
& [ \Pi (\xi ) , \phi (\eta ) ]\; = \; \frac{-i}{2} \; \left( \frac{}{} (\xi , \, \eta ) + (\eta , \xi ) \right) , \qquad  \xi , \eta \in \eltwo. \label{PiPhi}
\end{align}
By (\ref{PiHb}) and (\ref{PiPhi}), we have
\begin{equation} 
U (\Lambda )^{-1} 
 H (\Lambda  )  U (\Lambda ) 
   \kuuhaku = \kuuhaku  H_{0} (\Lambda )  \; 
    + K (\Lambda ) 
\end{equation}
where
\begin{equation}
K (\Lambda ) \; = \; 
U (\Lambda )^{-1}
  \left( \Hpar  \tens I \right)  U (\Lambda ) \;  -\;  \Hpar \tens I \; 
 \;  + \; 
  V_{\effect} (\mbf{x}_{1}, \cdots ,   \mbf{x}_{N}) \label{KLambda} .
\end{equation}

$\quad$ \\
$\; $ \\
Now we apply the general theory on scaling limits investigated in \cite{Ar90}.
Let us set the total Hilbert space  by $ \ms{Z} = \ms{X} \tens   \ms{Y}$.
Let $A$ and $B$ be non-negative self-adjoint
 operators on $\ms{X}$ and $\ms{Y}$, respectively.
Here we assume that ker $B \ne \{ 0 \} $. 
We consider a family of symmetric operators
 $\{ C(\Lambda ) \}_{\Lambda >0 } $ satisfying the 
 conditions :    
\begin{quote}
  (\textbf{S.1}) For all $\epsilon >0 $ there exists a constant 
 $\Lambda (\epsilon ) >0 $ such that for all 
$\Lambda > \Lambda (\epsilon )$, \\
 $\ms{D} (A \tens I) \cap 
 \ms{D} (I \tens B)  \subset \ms{D} (C(\Lambda ))$, 
 and there exists $ b(\epsilon ) \geq 0 $ such that
 \[
 \| C (\Lambda ) \Phi  \| \leq 
  \epsilon \| ( A \tens I  + \Lambda I \tens B ) \Phi \| 
   + b (\epsilon ) \| \Phi \| .
 \]
 \end{quote}
\begin{quote}
(\textbf{S.2}) There exists a symmetric operator $C$ on $\ms{Z}$ 
 such that   
 $ \ms{D} \tens \text{ker } B   \subset \ms{D} (C) $ 
 and for all $z \in \mbf{C} \backslash \mbf{R} )$, 
\[
s-\lim_{\Lambda \to \infty} 
C (\Lambda ) (A\tens I + \Lambda I \tens B -z) 
= C (A-z)^{-1} \tens P_{B},
\]  
where $P_{B} $ is the orthogonal
 projection from $\ms{Y}$ onto ker $B$. 
\end{quote}
\textbf{Proposition A} (\cite{Ar90} ; Theorem 2.1)  $\; $ Assume (\textbf{S.1}) and (\textbf{S.2}). Then (i)-(iii) follows.\\
$\quad$(i) There exists $\Lambda_{0} \geq 0$ such that for all 
 $\Lambda > \Lambda_{0} $, 
 \[
 X( \Lambda ) = A \tens I + \Lambda I \tens B 
 + C (\Lambda )
 \]is self-adjoint on $\ms{D}(A \tens I ) \cap \ms{D} 
(  I \tens B ) $ and uniformly bounded from below for 
 $\Lambda$, furthermore  $X(\Lambda )$ is essentially self-adjoint on any core of $A \tens I + I \tens B $. \\ 
$\quad$ (ii) Let $X= A \tens I + ( I \tens P_{B} )  C  ( I \tens P_{B} ) $.
 Then $X$ is self-adjoint on $\ms{D} (A \tens I )$ and  bounded from below, and essentially self-adjoint on any core of $A \tens I$. \\
$\quad$(iii) Let $z \in \bigcap_{\Lambda \geq \Lambda_{0}} 
 \rho (X(\Lambda )) \cap \rho (X) $, where $\rho (\ms{O}) $ denotes the resolvent set of an operator $\ms{O}$.  Then 
\[
s-\lim_{\Lambda \to \infty} (X(\Lambda ) -z)^{-1}
 = (X-z)^{-1} (I \tens P_{B}) . 
\]

$\quad $ \\
$\quad$ \\
$\quad $Now we consider $H(\Lambda) $ again. What we have to prove is that  
 $H (\Lambda )$ satisfies the condition \textbf{(S.1)} and \textbf{(S.2)} by applying
 $H_{0}(\Lambda )$ to $A \tens I + \Lambda \; I \tens B $ and $K (\Lambda) $ to $C (\Lambda)$. 
 First let us  consider  the term $U (\Lambda )^{-1} \left( \Hpar \tens I \right) U (\Lambda ) $ in (\ref{KLambda}).
 Let us set $\hat{\mbf{p}} \; = \; 
 ( \hat{p}^1 ,  \cdots , \hat{p}^d ) \; = \;
 (  -i \frac{\partial}{ \partial x^1} , \cdots , -i \frac{\partial}{ \partial x^d}  ) $.  
Then by the spectral decomposition theorem, 
\begin{equation}
U (\Lambda )^{-1}
  \left( \Hpar \tens I \right) U (\Lambda ) \; = \; 
  \sum_{j=1}^{N}
  \sqrt{\left(  \frac{}{} U(\Lambda )^{-1} ( \mbf{p}_{j} \tens I ) U (\Lambda ) \right)^2 \; \; + \; \; M^2 } ,
 \end{equation} 
follows. 
We see that 
\begin{equation}
[ \Pi (f_{\mbf{x}}) , \, \hat{p}^{\nu}  ] \; = \; i \Pi (\partial_{x^\nu} f_{\mbf{x}} ) .
\end{equation}
Then by (\textbf{A.4}), it follows that for $ \Psi \in \ms{D}_{0} $, 
\[
\left(  \frac{}{} U(\Lambda )^{-1} ( \mbf{p}_{j} \tens I ) U (\Lambda ) \right)^2  \Psi \; \; = \; \; 
   \left( \sum_{\nu=1 }^{d} 
\left( \phattwo{j}{\nu}  \tens I  + \left( \frac{\kappa}{\Lambda} \right)
 \Pi ( \frac{  \partial_{x_{j}^\nu }f_{\mbf{x}_{j}} }{ \omega} )
\right)^2  \; + \;  M^2 \right) \Psi .
\]
   Then   we have
 \begin{equation}
\left(  \frac{}{} U(\Lambda )^{-1} ( \mbf{p}_{j} \tens I ) U (\Lambda ) \right)^2 
 \; = \; \overline{
   \left( \sum_{\nu=1 }^{d} 
\left( \phattwo{j}{\nu}  \tens I  + \left( \frac{\kappa}{\Lambda} \right)
 \Pi ( \frac{  \partial_{x_{j}^\nu }f_{\mbf{x}_{j}} }{ \omega} )
\right)^2  \; + \;  M^2 \right)_{\upharpoonright \ms{D}_{0}} }  ,
\end{equation}
where $\overline{Y}$ denotes the closure of the operator $Y$.
 Here we abbreviate  as
\begin{align*}
& \Pi ( \frac{\mbf{\nabla} f_{\mbf{x}} }{ \omega} )
 \cdot ( \hat{\mbf{p}} \tens I ) \; 
= 
\; \sum_{\nu =1}^{d} 
 \Pi ( \frac{ \partial_{x^\nu} f_{\mbf{x}}}{ \omega} )
( \phatone{\nu}  \tens I ),  \\
&  \Pi ( \frac{\mbf{\nabla} f_{\mbf{x}}}{ \omega} ) \cdot  
\Pi ( 
\frac{\mbf{\nabla} f_{\mbf{x}}}{ \omega} )
\; = \; 
\sum_{\nu = 1 }^{d} \Pi ( \frac{ \partial_{x^\nu} f_{\mbf{x}}}{ \omega} )
 \phi ( \frac{ \partial_{x^\nu} f_{\mbf{x}}}{ \omega} ) . 
\end{align*}
Then we see that
\begin{equation}
U (\Lambda )^{-1}
 \left( \Hpar \tens I \right) U (\Lambda ) \; = \; 
  \sum_{j=1}^{N} \left(
   \; \overline{ \;  
\left( \frac{}{} - \triangle_{j}  \tens I \; + \; Q_{j} (\Lambda )  \; + M^2 \right)_ {\upharpoonright \ms{D}_{0}}} \right)^{1/2}  ,  \label{4/12.1}
\end{equation}
where
\[
Q_{j} (\Lambda ) \; = \; 
\left( \frac{\kappa }{ \Lambda } \right)
 \left(
2 \Pi ( \frac{\mbf{\nabla} f_{\mbf{x}_j } }{ \omega} )
 \cdot ( \hat{\mbf{p}_j} \tens I ) \; - \; i \Pi (\frac{\triangle f_{\mbf{x}} }{\omega} )
 \right)
 + 
\left( \frac{\kappa}{\Lambda} \right)^2 \Pi ( \frac{\mbf{\nabla} f_{\mbf{x}_j}}{ \omega} ) \cdot  
\Pi ( 
\frac{\mbf{\nabla} f_{\mbf{x}_j}}{ \omega} )
 .
\]

\begin{proposition} \label{mainprop}
Assume \textbf{(A.1)}-\textbf{(A.4)}. Then for $\epsilon > 0$, 
there exists $\Lambda (\epsilon) \geq 0 $ such that 
 for all  $  \Lambda > \Lambda (\epsilon) $,
\begin{equation}
\| U (\Lambda )^{-1}  \left( \Hpar  \tens I \right) U (\Lambda )  \Psi
 - \left( \Hpar  \tens I \right ) \Psi  \| \; 
\leq \;  \epsilon \| H_{0} (\Lambda ) \Psi \| \,  +
  \,  b (\epsilon ) \| \Psi \| 
  \label{6/17.2}
\end{equation}  
 where $ b (\epsilon )$ is a constant independent of 
 $ \Lambda \geq \Lambda (\epsilon) $.
\end{proposition}

$\quad$ \\
Before proving Proposition \ref{mainprop}, we show the following lemma.
\begin{lemma} \label{mlemma}
For $\lambda > 0$ and $\delta \in (0, \frac{1}{10})$, there exists $M_{\nu} (\delta )$, $\nu =1, \cdots , d$, such that  
 \begin{equation}
\| \hat{p}^{\nu} ( - \triangle + M^2 
+ \lambda )^{-1} ( \sqrt{-\triangle + M^2 } + 1 )^{-1/2} \| \leq \; 
 \frac{1}{\lambda^{\frac{1}{2} + \delta}} M_{\nu} (\delta ) .   \label{mainlemma}
  \end{equation}
\end{lemma}
(\textbf{Proof})
For $\mbf{p} = (p^{1} , \cdots , p^{d} ) \in \mbf{R}^d$, $\nu = 1, \cdots , d$, 
we see that 
\begin{equation}
|\; p^\nu ( \mbf{p}^2 + M^2 + \lambda  )^{-1}  (\sqrt{\mbf{p}^2 + M^2 } + 1 )^{-1/2} \;  |
 =\frac{1}{\lambda^{\frac{1}{2}+ \delta}}  \;  |  \lambda^{\frac{1}{2}+ \delta}p^\nu | 
 \;  ( \mbf{p}^2 + M^2 + \lambda  )^{-1}  (\sqrt{\mbf{p}^2 + M^2 } +1 )^{-1/2} . \notag
\end{equation}
We shall show that 
\begin{equation}
\sup_{\lambda > 0 , \,\mbf{p} \in \Rd } \; 
|  \lambda^{\frac{1}{2}+ \delta}p^\nu | \;
 \;  ( \mbf{p}^2 + M^2 + \lambda  )^{-1}  (\sqrt{\mbf{p}^2 + M^2 } +1 )^{-1/2} \; < \; \infty , \label{3/30.0}
\end{equation}
and hence  (\ref{mainlemma}) follows from the spectral decomposition theorem. 
The Young's inequality shows that 
for $q>1$ and   $\tilde{q} > 1 $ satisfying $  \frac{1}{q} + \frac{1}{\tilde{q}} = 1 $,  
\begin{equation}
\lambda^{\frac{1}{2} + \delta } |p^\nu | \; \leq 
\;  \frac{1}{q} \lambda^{ ( \frac{1}{2} + \delta ) q }
  \; + \;  \frac{1}{\tilde{q}} |p^\nu |^{\tilde{q}}  
\end{equation}
follows. Let us take  $q = ( \frac{1}{2} + \delta  )^{-1}$ for $\delta \in (0, \frac{1}{10})$, and hence 
$\; \tilde{q} = ( \frac{1}{2} - \delta  )^{-1} $. Then we have 
\begin{equation}
\lambda^{\frac{1}{2} + \delta } |p^\nu|
\; \leq \; 
 ( \frac{1}{2} + \delta  ) \lambda \; + \; 
  ( \frac{1}{2} - \delta  ) |p^\nu |^{( \frac{1}{2} - \delta  )^{-1}}. \label{3/30.2}
\end{equation}
Note that 
\begin{equation}
\sup_{\lambda >0 , \mbf{p} \in \Rd } 
 \; \; \lambda  \;  ( \mbf{p}^2 + M^2 + \lambda  )^{-1}  (\sqrt{\mbf{p}^2 + M^2 } +1 )^{-1/2} \; < \; \infty .
 \label{3/30.3}
\end{equation}
Since  $ 0< \delta < \frac{1}{10}$, we see that 
$   ( \frac{1}{2} - \delta  )^{-1} \; < \frac{5}{2} $, and hence
\begin{equation}
\sup_{\mbf{p} > \Rd}  |p^\nu |^{( \frac{1}{2} - \delta  )^{-1}} 
 (\mbf{p}^2  + M^2   )^{-1} (\sqrt{\mbf{p}^2 +M^2} +1 )^{-1/2}
< \infty . 
\end{equation}
Then we have 
\begin{align}
& \sup_{\lambda >0 , \mbf{p} \in \Rd }  |p^\nu |^{( \frac{1}{2} - \delta  )^{-1}} 
 (\mbf{p}^2  + M^2  + \lambda )^{-1} (\sqrt{p^2 +M^2} +1 )^{-1/2}  \notag \\
  &\qquad \qquad \qquad  \leq  \sup_{\mbf{p} \in \Rd }  |p^\nu |^{( \frac{1}{2} - \delta  )^{-1}} 
 (\mbf{p}^2  + M^2   )^{-1} (\sqrt{\mbf{p}^2 +M^2} +1 )^{-1/2}
< \infty .  \label{3/30.4}
\end{align}
By (\ref{3/30.2}), (\ref{3/30.3}) and (\ref{3/30.4}), we obtain (\ref{3/30.0}). $\blacksquare$

$\quad$ \\
\textbf{{\large (Proof of Proposition \ref{mainprop})}}  \\
It follows that for a nonnegative self-adjoint operator $S$,
\begin{equation}
\sqrt{S} \Phi \; = \; 
 \frac{1}{\pi} \int_{0}^{\infty}
 \frac{1}{\sqrt{\lambda}} (S + \lambda  )^{-1} \,  S \,  \Phi
 \;  d \lambda
 , \qquad \qquad \Phi \in \ms{D} (S) .
\end{equation}
Let 
\begin{align*}
&A_{j} (\Lambda ) \; =  \; \overline{\left( \frac{}{} - \triangle_{j} \tens I \; +
  Q_{j} (\Lambda )\; +  \; M^2 \right)_{\upharpoonright \ms{D}_{0}}} ,  \\  
&B_{j} \; = \;  - \triangle_{j} \tens I \; + \; M^2  .
\end{align*}
Then we  have for $ \Psi  \in \ms{D}_{0}  $, 
\begin{align}
\left( \frac{}{} U (\Lambda )^{-1} \left( \Hpar \tens I \right)  U (\Lambda )   - \Hpar \tens I  \right) \Psi  
&= \sum_{j=1}^{N} \frac{1}{\pi} \int_{0}^{\infty}
 \frac{1}{\sqrt{\lambda}} 
\left\{  (A_{j} (\Lambda ) + \lambda  )^{-1} A_{j} (\Lambda )\;
- \; (B_{j} + \lambda )^{-1} B_{j} \right\} \,  \Psi 
 \;  d \lambda  \notag \\ 
 &=  \sum_{j=1}^{N} \frac{1}{\pi} \int_{0}^{\infty}
 \sqrt{\lambda} 
  (A_{j} (\Lambda ) + \lambda  )^{-1}
 ( A_{j} (\Lambda )  
- \; B_{j} )  (B_{j} + \lambda )^{-1}  \,  \Psi 
 \;  d \lambda  \notag \\ 
&  =\sum_{j=1}^{N} \frac{1}{\pi} \int_{0}^{\infty}
\sqrt{\lambda} 
  (A_{j} (\Lambda ) + \lambda  )^{-1}
 \; Q_{j}(\Lambda ) \; (B_{j} + \lambda )^{-1}  \,  \Psi 
 \;  d \lambda .
\end{align}
By (\ref{4/12.1}) and the spectral decomposition theorem, $ \|  (A_{j} (\Lambda ) + \lambda  )^{-1} \|  \leq 
 \frac{1}{\lambda + M^2} $, $\; \lambda >0 $ follows, and then  we have
\begin{equation}
\| \left( \frac{}{} U (\Lambda )^{-1}  \left( \Hpar \tens I \right)   U (\Lambda )   -  \Hpar  \tens I \right) \Psi \| \leq 
 \sum_{j=1}^{N} \frac{1}{\pi} \int_{0}^{\infty}
\frac{\sqrt{\lambda}}{\lambda + M^2} 
   \; \|  Q_{j}(\Lambda ) \; (B_{j} + \lambda )^{-1}  \,  \Psi  \|
 \;  d \lambda . \label{hoshi0}
\end{equation}
We see that
\begin{align}
  \| Q_{j}(\Lambda ) \; (B_{j} + \lambda )^{-1}  \,  \Psi \| 
   &\leq  \left( \frac{\kappa}{\Lambda} \right) \left(  \|  \Pi ( \frac{\mbf{\nabla} f_{\mbf{x}_j } }{ \omega} )
 \cdot ( \hat{\mbf{p}_j} \tens I ) (B_{j} + \lambda )^{-1} \Psi \| + 
 \| \Pi (\frac{\triangle f_{\mbf{x}_{j}} }{\omega} ) (B_{j} + \lambda )^{-1} \Psi \| \right) \notag \\
 & \quad +  \left( \frac{\kappa}{\Lambda}\right)^2  \| \Pi ( \frac{\mbf{\nabla} f_{\mbf{x}_j}}{ \omega} ) \cdot  \Pi ( \frac{\mbf{\nabla} f_{\mbf{x}_j}}{ \omega} )
 (B_{j} + \lambda )^{-1}  \Psi \|  . \label{3/31.3}
\end{align}
Note that 
\begin{align}
& \|  \Pi ( \frac{\mbf{\nabla} f_{\mbf{x}_j } }{ \omega} )
 \cdot ( \hat{\mbf{p}_j} \tens I ) (B_{j} + \lambda )^{-1} \Psi \|  \notag \\
 & \leq \sum_{\nu } 
 \| \Pi ( \frac{ \partial_{x_j^\nu} f_{\mbf{x}_j } }{ \omega} ) ( I \tens \Hb + 1)^{-1/2} \| \; 
 \|   ( \hat{\mbf{p}_j}^\nu \tens I ) (B_{j} + \lambda )^{-1} ( \Hpar \tens I + 1 )^{-1/2} \| \; 
\| ( \Hpar \tens I  + 1 )^{1/2}( I \tens \Hb  + 1)^{1/2}  \Psi \| \label{4/6.1} 
\end{align}
Here we   used the boundness (\ref{bounda}) and (\ref{boundad}). Applying the Lemma \ref{mlemma} to $ \|   ( \hat{\mbf{p}_j}^\nu \tens I ) (B_{j} + \lambda )^{-1} ( \Hpar \tens I   + 1 )^{-1/2} \| $ in (\ref{4/6.1}), 
   it is seen that for $\delta \in (0, \frac{1}{10})$, there exist $\alpha_{j} (\delta) \geq 0  $  such that 
\begin{equation}
 \|  \Pi ( \frac{\mbf{\nabla} f_{\mbf{x}_j } }{ \omega} )
 \cdot ( \hat{\mbf{p}_j} \tens I ) (B_{j} + \lambda )^{-1} \Psi \| 
 \leq \frac{\alpha_{j}(\delta )}{\lambda^{\frac{1}{2} + \delta }}  
\| ( \Hpar \tens I  + 1 )^{1/2}( I \tens \Hb + 1)^{1/2}  \Psi \| \notag   ,
\end{equation}
and hence we have
\begin{equation}
 \|  \Pi ( \frac{\mbf{\nabla} f_{\mbf{x}_j } }{ \omega} )
 \cdot ( \hat{\mbf{p}_j} \tens I ) (B_{j} + \lambda )^{-1} \Psi \| \leq \frac{\alpha_{j} (\delta )}{\lambda^{\frac{1}{2} + \delta }}  
 \left( \frac{}{}  \|  \Hpar   \Psi   \| + \| \Hb  \Psi \| + \| \Psi \| \right) .  \label{hoshi1}
\end{equation}
Since $ (\|  \Hpar  \tens I \Psi   \| + \|  I \tens \Hb  \Psi \|)^2 \leq 2 \| H_0 (\Lambda) \|^2 $,
we have 
\begin{equation}
\int_{0}^{\infty }  \frac{\sqrt{\lambda}}{\lambda + M^2}  \|  \Pi ( \frac{\mbf{\nabla} f_{\mbf{x}_j } }{ \omega} )
 \cdot ( \hat{\mbf{p}_j} \tens I ) (B_{j} + \lambda )^{-1} \Psi \| d \lambda  
 \leq \alpha_{j} (\delta )\left( \int_{0}^{\infty } \frac{1}{ (\lambda +M^2 )\lambda^{ \delta }} d \lambda \right)
 \left( \frac{}{} \sqrt{2} \|  H_{0} (\Lambda ) \Psi   \|  + \| \Psi \| \right) . \label{3/31.a}
\end{equation}
By   $\| \Pi (\frac{\triangle f_{\mbf{x}_{j}} }{\omega} ) (I \tens \Hb + 1 )^{-1/2} \| < \infty \; $ 
 and $ \;  \|  (B_{j} + \lambda )^{-1} \| \leq  \frac{1}{\lambda + M^2 }  $, we have 
\begin{align}
\| \Pi (\frac{\triangle f_{\mbf{x}_{j}} }{\omega} ) (B_{j} + \lambda )^{-1} \Psi \| 
& \leq  \| \Pi (\frac{\triangle f_{\mbf{x}_{j}} }{\omega} ) (I \tens \Hb + 1 )^{-1/2} \| \; 
 \|  (B_{j} + \lambda )^{-1} \| \;   \|   (I \tens \Hb + 1 )^{1/2}  \Psi \|   \notag \\
& \leq \frac{1}{\lambda + M^2 }  \| \Pi (\frac{\triangle f_{\mbf{x}_{j}} }{\omega} ) (I \tens \Hb + 1 )^{-1/2} \|  \; 
\|   (I \tens \Hb + 1 )^{1/2}  \Psi \|  .\label{hoshi2}
\end{align}
Then  by  $\;  \|   ( I \tens \Hb + 1 )^{1/2}  \Psi \|  \leq \| H_{0} (\Lambda ) \Psi \| + \| \Psi \|$,  we have
\begin{align}
 \int_{0}^{\infty } \frac{\sqrt{\lambda}}{\lambda + M^2}  \| \Pi (\frac{\triangle f_{\mbf{x}_{j}} }{\omega} ) & (B_{j} + \lambda )^{-1} \Psi \|  d \lambda \notag \\
 & \leq   \| \Pi (\frac{\triangle f_{\mbf{x}_{j}} }{\omega} ) (I \tens \Hb + 1 )^{-1/2} \|\left( \int_{0}^{\infty } 
\frac{\sqrt{\lambda}}{  (\lambda + M^2 )^{2} } d \lambda \right)
 \left( \frac{}{}  \|  H_{0} (\Lambda ) \Psi   \|  + \| \Psi \| \right) . \label{3/31.b}
\end{align} 
In addition we  also see that
\begin{align}
 \| \Pi ( \frac{\mbf{\nabla} f_{\mbf{x}_j}}{ \omega} ) \cdot  
\Pi ( \frac{\mbf{\nabla} f_{\mbf{x}_j}}{ \omega} )
 (B_{j} + \lambda )^{-1}  \Psi \|  
 & \leq \| \Pi ( \frac{\mbf{\nabla} f_{\mbf{x}_j}}{ \omega} ) \cdot  
\Pi ( \frac{\mbf{\nabla} f_{\mbf{x}_j}}{ \omega} ) ( I \tens \Hb + 1 )^{-1} \| \; 
 \|  (B_{j} + \lambda )^{-1} \| \; \|  ( I \tens \Hb + 1 )\Psi \|  \notag \\
 & \leq \frac{1}{\lambda + M^2 } \| \Pi ( \frac{\mbf{\nabla} f_{\mbf{x}_j}}{ \omega} ) \cdot  
\Pi ( \frac{\mbf{\nabla} f_{\mbf{x}_j}}{ \omega} ) ( I \tens \Hb + 1 )^{-1} \|  \; \|  ( I \tens \Hb +1)\Psi  \| . 
\label{hoshi3}
\end{align}
Since $\| \Pi (\xi ) \Pi (\eta ) (\Hb +1 )^{-1} \| <\infty $  for $\xi , \eta \in {\ms{D} (\omega ) }$, we obtain
\begin{align}
 \int_{0}^{\infty }  & \frac{\sqrt{\lambda}}{\lambda + M^2}  \| \Pi ( \frac{\mbf{\nabla} f_{\mbf{x}_j}}{ \omega} ) \cdot  
\Pi ( \frac{\mbf{\nabla} f_{\mbf{x}_j}}{ \omega} )
 (B_{j} + \lambda )^{-1}  \Psi \|  d \lambda  \notag \\
 & \quad  
 \leq    \| \Pi ( \frac{\mbf{\nabla} f_{\mbf{x}_j}}{ \omega} ) \cdot  
\Pi ( \frac{\mbf{\nabla} f_{\mbf{x}_j}}{ \omega} ) ( \Hb + 1 )^{-1} \| 
|\left( \int_{0}^{\infty } \frac{\sqrt{\lambda}}{(\lambda + M^2 )^2  } d \lambda \right)
 \left( \frac{}{}  \|  H_{0} (\Lambda ) \Psi   \|  + \| \Psi \| \right) . \label{3/31.c}
\end{align} 
Then from  (\ref{3/31.a}),(\ref{3/31.b}), (\ref{3/31.c}) and (\ref{3/31.3}), the proposition follows. 
$\blacksquare$

%%%%%%%%%%%%%%%%%%%%%%%%%%%%%%%%%%%%%%%%%%%%%%%%%%%%%%%%%%%%%%%%%%%%%%%%%%%%%%%%%%%%%%%%%%%%%%%%%
\begin{proposition} \label{3/31.1}
  Assume \textbf{(A.1)} - \textbf{(A.4)}. \\
\textbf{(1)}$\; $ For $\epsilon > 0$, 
there exists $\Lambda (\epsilon) \geq 0 $ such that 
 for all  $  \Lambda > \Lambda (\epsilon) $,
\begin{equation}
\qquad \| K(\Lambda )  \Psi \| \; 
\leq \;  \epsilon \| H_{0} (\Lambda ) \Psi \| \,  +
  \,  \nu (\epsilon ) \| \Psi \|  , \qquad \qquad \Psi \in \ms{D}_{0} , 
  \label{6/17.2}
\end{equation}  
 holds, where $ \nu (\epsilon )$ is a constant independent of 
 $ \Lambda \geq \Lambda (\epsilon) $. \\ 
\textbf{(2)}$\; $ 
Then for all $ z  \in  \mbf{C} \backslash  \mbf{R} $, it follows that 
\begin{equation}
 s-\lim_{\Lambda \to \infty }
K(\Lambda ) 
\left( H_{0} (\Lambda ) -z \right)^{-1}
 \;   =   \;  V_{\effect} \; ( \;  \Hpar \; - \; z \;    )^{-1} \tens P_{\Omega_{\bos}}  ,
 \label{6/17.5}
 \end{equation}
 \end{proposition}
(\textbf{Proof}) \\
\textbf{(1)} By the condition (\textbf{A.4}), $V_{\effect} $ is bounded. 
 Then \textbf{(1)} follows from Proposition \ref{mainprop}.
$\;$ \\
\textbf{(2)}
It is seen that
\[
K(\Lambda ) 
\left( H_{0} (\Lambda ) - \frac{}{}z  \right)^{-1} \;
 \;  = \; K(\Lambda ) 
\left(   \Hpar  - \frac{}{} z  \right)^{-1}
 \tens \POmegab     
+ K(\Lambda ) 
\left( H_{0} (\Lambda ) - \frac{}{}z  \right)^{-1}  
 \left(  I \tens  (1 - \POmegab) \frac{}{} \right) .
\]
By Proposition \ref{mainprop},  we have 
\begin{equation}
s-\lim_{\Lambda \to \infty }
 K(\Lambda )  \left( \frac{}{}
\left( \frac{}{}  \Hpar -  z  \right)^{-1}
 \tens \POmegab  \right)  \Psi  \; = V_{\effect}
 \left( \frac{}{}
\left(  \Hpar - \frac{}{} z  \right)^{-1}
 \tens \POmegab \right) \Psi  .   \label{6/17.3}
\end{equation}
By  (\ref{6/17.2}), we see that 
 for $\epsilon >0 $ there exists $ \Lambda (\epsilon ) \geq 0 $ such that for all $ \Lambda > \Lambda (\epsilon )$  
\begin{equation}
\|  K(\Lambda ) 
\left( H_{0} (\Lambda ) - \frac{}{}z  \right)^{-1}  \Phi \| 
\leq  \epsilon \| \Phi \| + 
( \epsilon |z|  + \nu (\epsilon  ) ) \|  \left( H_{0} (\Lambda ) - \frac{}{}z  \right)^{-1}  \Phi \| , \qquad\; \Phi \in \ms{H} . \notag
\end{equation}
Note that
$  \lim_{\Lambda \to \infty} \; \left\|  \left( H_{0} (\Lambda ) - \frac{}{}z  \right)^{-1}  
  \left(  I \tens  (1 - \POmegab ) \frac{}{} \right)  \Psi \right\|  \; = \; 0 $, 
and hence we obtain 
\begin{equation}
\lim_{\Lambda \to \infty } \left\| 
 K(\Lambda ) 
\left( H_{0} (\Lambda ) - \frac{}{}z  \right)^{-1}  
 \left(  I \tens   (1 - \POmegab) \frac{}{} \right)   \Psi  \right\|
  \; = \; 0  . \label{6/17.4}
\end{equation}
By (\ref{6/17.3}) and (\ref{6/17.4}), 
 we obtain (\ref{6/17.5}).  
$\blacksquare $

$\quad$ \\
\textbf{{\large (Proof of Theorem \ref{MainTheorem})}} \\
By Proposition \ref{3/31.1}, it is shown that 
$H (\Lambda )$ satisfies the condition \textbf{(S.1)} and \textbf{(S.2)} by applying
 $H_{0}(\Lambda )$ to $A \tens I + \Lambda \; I \tens B $ and $K (\Lambda) $ to $C (\Lambda)$. 
Hence by the Proposition A, we have for $z \in \mbf{C} \backslash \mbf{R}$, 
\[
s-\lim_{\Lambda \to \infty} \left( H (\Lambda) - z \frac{}{} \right)^{-1}  \; 
= \;  \lim_{\Lambda \to \infty } U (\Lambda )\left( H_{0} (\Lambda ) - z\frac{}{} \right)^{-1}    U (\Lambda )^{-1}
 \; = \left( \Hpar + V_{\eff} -z \frac{}{} \right)^{-1} \tens \POmegab .
 \]
Thus the proof is completed. $\blacksquare$ \\

$\quad$ \\
{\Large Acknowledgments} \\
It is pleasure to thank Professor  Fumio Hiroshima for his advice and comments.


\begin{thebibliography}{99}
\footnotesize{
 \bibitem{Ar90}
A. Arai, Asymptotic analysis and its application to the nonrelativistic
  limit of the Pauli-Fierz and a spin-boson model,
  {\it J. Math. Phys.} \textbf{32} (1990),  2653-2663.
 \bibitem{Da79}
E. B. Davies, Particle-boson interactions and the weak coupling limit,
{\it  J. Math. Phys.} \textbf{20} (1979),  345-351.
 \bibitem{Da84}
 I. Daubechies, One electron molecules with relativistic kinetic energy : Properties of the discrete spectrum, 
 \textit{Commu. Math. Phys.}, \textbf{94} (1984), 523-535.
 \bibitem{De90}
J. Derezi\'nski, The Mourre estimate for dispersive $N$-body Schr\"odinger operators, 
 \textit{Trans. Amer. Math. Soc.} \textbf{317} (1990), 773-798.
 \bibitem{Ge91} 
C. G\'erard, The mourre estimate for regular dispersive systems,
 \textit{Ann. Inst. H. Poincar\'e Phys. Th\'eor.} \textbf{54} (1991), 59-88. 
 \bibitem{He77}
I. Herbst, Spectral theory of the operator $(p^2 + m^2 )^{1/2} - Ze^2 / r$, 
\textit{Commun. Math. Phys.}, \textbf{53} (1977), 285-294.
 \bibitem{Hiro}
F.  Hiroshima,  Analysis of ground states of atoms interacting with a quantized radiation field,
\textit{ Topics in the theory of  Shr\"{o}dinger operators },
 (H.Araki and H.Ezawa eds. ) World Scientific, 2004,
 145-273.
 \bibitem{Hi93}
 F.  Hiroshima,  Scaling limit of a model of quantum electrodynamics,
 {\it J. Math. Phys.} \textbf{34} (1993),  4478-4518.
 \bibitem{Hi97}
 F.  Hiroshima,  Scaling limit of a model of quantum electrodynamics with many nonrelativistic particles,
 {\it Rev. Math. Phys.} \textbf{9} (1997),  201-225.
 \bibitem{Hi98}
 F.  Hiroshima,  Weak coupling limit with a   removal of an ultraviolet   cutoff  for a Hamiltonian of particles interacting with a massive scalar field,
 {\it Inf. Dim. Ana. Quantum Prob. Rel. Top.} \textbf{1} (1998),  407-423.
 \bibitem{Hi99}
 F.  Hiroshima,  Weak coupling limit and a removing  ultraviolet cutoff for a Hamiltonian of particles interacting with a   quantized scalar field,
 {\it  J. Math. Phys. } \textbf{40} (1999),
 1215-1236.
 \bibitem{Hi02}
 F.  Hiroshima,  Observable effects and parametrized scaling limits of a model in non-relativistic quantum electrodynamics,
{\it  J. Math. Phys.} \textbf{43} (2002),  1755-1795.
\bibitem{HiSp01}
 F.  Hiroshima  and H.Spohn, Enhanced binding through  coupling to a quantum field, {\it  Ann. Henri. Poincar\'{e}} \textbf{2}
  (2001),  1150-1187.
 \bibitem{HiSa}
F.  Hiroshima  and I.Sasaki, On the ionization of the semi-relativistic Pauli-Fierz model for a single particle, 
 (arxiv : 1003.1661v4).
 \bibitem{LiLo}
E. Lieb and M. Loss, \textit{Analysis} (second edition), Amer. Math. Soc. 2001.
\bibitem{Oh09}
 A. Ohkubo, Scaling limit for the Derezin\'{n}ski-G\'{e}rard Model, to appear in Hokkaido Math. J.
  \bibitem{Su07-1}
 A. Suzuki, Scaling limits for a general class of quantum field models and its applications  to nuclear physics and condensed matter physics, {\it Inf. Dim. Ana. Quantum Prob. Rel. Top.}  \textbf{10} (2007),  43-65.
\bibitem{Su07-2}
 A. Suzuki, Scaling limits for a generalization of the Nelson model and its application to nuclear physics, {\it Rev. Math. Phys.}
  \textbf{19} (2007),  131-155.
  \bibitem{SSS10}
 J. P. Solovej, T. \O. S{\o}rensen, and W. L. Spitzer, Relativistic Scott correction for atoms and molecules, 
{\it Comm. Pure Appl. Math.} \textbf{63} (2010), 39-118.  
 \bibitem{Ta09}
 T. Takaesu, On the scaling limit of quantum electrodynamics with spatial cutoffs. (arxiv : 0908.2080v1)
 \bibitem{Um95}
 T. Umeda, Radiation conditions and resolvent estimates for relativistic Schr\"odinger operators, 
  \textit{Ann. Inst. H. Poincar\'e Phys. Th\'eor.} \textbf{63} (1995), 277-296. 
\bibitem{We74}
 R. Weder,  Spectral properties of one-body relativistic spin-zero hamiltonians,  
\textit{Ann. Inst. H. Poincar\'{e}, Sect. A}, \textbf{20} (1974), 211-220.
} 
 
\end{thebibliography}
\end{document}